\documentclass[12pt,twoside]{article}
\usepackage[latin1]{inputenc}
\usepackage{amsmath}
\usepackage{amssymb,amsfonts}
\usepackage{graphicx}
\usepackage{times,amssymb,amscd}
\usepackage[USenglish]{babel}
\usepackage{verbatim}
\usepackage{enumerate}

\newcommand{\bH}{\mathbf{H}}

\newcommand{\cD}{\mathcal{D}}

\newcommand{\cB}{\mathcal{B}}

\begin{document}
\pagestyle{myheadings}
\markboth{\centerline{J.~Szirmai}}
{Horoball packings and their densities...}
\title
{Horoball packings and their densities by generalized simplicial density function in the hyperbolic space
\footnote{Mathematics Subject Classification 2010: 52C17, 52C22, 52B15. \newline
Key words and phrases: Hyperbolic geometry, horoball packings, optimal simplicial density.}}

\author{\medbreak \medbreak {\normalsize{}} \\
\normalsize Jen\H{o}  Szirmai \\
\normalsize  Budapest University of Technology and Economics\\
\normalsize Institute of Mathematics, Department of Geometry \\
\normalsize H-1521 Budapest, Hungary \\
\normalsize Email: szirmai@math.bme.hu }
\date{\normalsize (\today)}
%%%%%%%%%%%%%%%%%%%%%%%%%%%%%%%%%%%%%%%%%%%%
%%%%%%%%%%%%%%%%%%%%%%%%%%%%%%%%%%%%%%%%%%%%
\maketitle
%%%%%%%%%%%%%%%%%%%%%%%%%%%%%%%%%%%%%%%%%%%%
\begin{abstract}
The aim of this paper is to determine the locally densest horoball packing arrangements and their densities with
respect to fully asymptotic tetrahedra with at least one plane of symmetry in hyperbolic 3-space $\overline{\mathbf{H}}^3$ extended with its 
absolute figure, 
where the ideal centers of horoballs give rise to vertices of a fully asymptotic tetrahedron. 
We allow horoballs of different types at the various vertices. Moreover, we generalize the notion of the simplicial density function in the 
extended hyperbolic space $\overline{\mathbf{H}}^n, ~(n \ge 2)$, and prove
that, in this sense, {\it the well known B\"or\"oczky--Florian 
density upper bound for "congruent horoball" packings of $\overline{\mathbf{H}}^3$
does not remain valid to the fully asymptotic tetrahedra.} 

The density of this locally densest packing is $\approx 0.874994$, may be surprisingly larger than the B\"or\"oczky--Florian 
density upper bound $\approx 0.853276$ but our local ball arrangement seems not to have extension to the whole hyperbolic space.  
\end{abstract}

%%%%%%%%%%%%%%%%%%%%%%%%%%%%%%%%%%%%%%%%%%%
%%%%%%%%%%%%%%%%%%%%%%%%%%%%%%%%%%%%%%%%%%%

\newtheorem{theorem}{Theorem}[section]
\newtheorem{corollary}[theorem]{Corollary}
\newtheorem{lemma}[theorem]{Lemma}
\newtheorem{exmple}[theorem]{Example}
\newtheorem{defn}[theorem]{Definition}
\newtheorem{rmrk}[theorem]{Remark}
\newtheorem{proposition}[theorem]{Proposition}
\newenvironment{definition}{\begin{defn}\normalfont}{\end{defn}}
\newenvironment{remark}{\begin{rmrk}\normalfont}{\end{rmrk}}
\newenvironment{example}{\begin{exmple}\normalfont}{\end{exmple}}
\newenvironment{acknowledgement}{Acknowledgement}

%%%%%%%%%%%%%%%%%%%%%%%%%%%%%%%%%%%%%%%%%%%%%%%%%%%%%%%%%%%%%%%%%%%%

%==================================================================%
%                             the main article                               %
%==================================================================%

%%%%%%%%%%%%%%%%%%%%%%%%%%%%%%%%%%%%%%%%%%%%%%%%%%%%%%%%%%%%%%%%%%%%

\section{Basic notions}

\subsection{Local density function}
We summarize the most important definitions and results about
ball packings in $\overline{\mathbf{H}}^n$ $(n\ge2)$. For more details and proofs, we refer to \cite{Be}, \cite{B78}, \cite{K98}, \cite{KSz} and \cite{Ro64}.
There are different notions of packing density. For later purposes, the local density measure
is the best suited one.
In the $n$-dimensional $(n \ge 2)$ hyperbolic space there are 3-types of spheres: sphere, horosphere and hypersphere.
Now, we consider the horospheres and their bodies, the horoballs. A horoball packing $\cB_h=\{B_h\}$ of $\overline{\mathbf{H}}^n$ is an arrangement of non-overlapping 
horoballs ${B_h}$ in $\overline{\mathbf{H}}^n$. The notion of local density of the usual ball packing can be extended for horoball packings $\cB_h$ of $\overline{\mathbf{H}}^n$. Let
$B_h \in \cB_h$, and $P \in \overline{\mathbf{H}}^n$ an arbitrary point. Then, $\rho(P,B_h)$ is defined to be the length of the
unique perpendicular from $P$ to the horosphere $S_h$ bounding $B_h$, where again $\rho(P,B_h)$ 
is taken negative for $P \in B_h$. The Dirichlet--Voronoi cell (shortly D-V cell) $\cD(B_h)$ of $B_h$ in $\cB_h$ is defined to be
the convex body
\begin{equation}
\cD_h = D(\cB_h,B_h) := \{ P \in \mathbf{H}^n ~|~ \rho(P,B_h) \le \rho(P,B'_h), ~ \forall B'_h \in \cB_h \}. \tag{1.1}
\end{equation}
Since both, $B_h$ and $\cD_h$, are of infinite volume, the usual concept of local density has to be
modified. Let $Q \in \partial{\mathbf{H}}^n$ denote the base point (ideal center at the infinity) of $B_h$, and interpret $S_h$ as a Euclidean
$(n - 1)$-space. Let $B_{n-1}(R) \subset S_h$ be an $n-1$-ball with center $C \in S_h$. Then,
$Q \in \partial {\mathbf{H}^n}$ and $B_{n-1}(R)$ determine a convex cone $C_n(R) :=cone(B^Q_{n-1}(R)) \in \overline{\mathbf{H}}^n$ with
apex $Q$ consisting of all hyperbolic geodesics through $B_{n-1}(R)$ with limiting point $Q$. With
these preparations, the local density $\delta_n(B_h, \cB_h)$ of $B_h$ to $\cD_h$ is defined by 
\begin{equation}
\delta_n(\cB_h, B_h) :=\overline{\lim_{R \rightarrow \infty}} \frac{vol(B_h \cap C_n(R))} {vol(\cD_h \cap C_n(R))}, \tag{1.2}
\end{equation}
and this limes superior is independent of the choice of the center $C$ of $B_{n-1}(R)$. 

\subsection{Densest packings with horoballs of the same type}

{\it We have to change the notion of the ,,congruent" horoballs in a horoball packing to the horoballs of the "same type" because the horroballs are 
congruent in the hyperbolic space $\overline{\mathbf{H}}^n$. Two horoballs in a horoball packing are in the "same type" if and only if the 
local densities of the horoballs to the corresponding cell (e.g. D-V cell; or ideal simplex, later on) are equal.} 
{\it \bf If we assume that the ,,horoballs belong to the same type"}, then by analytical continuation,
the well known simplicial density function on $\overline{\mathbf{H}}^n$ can be extended from $n$-balls of radius $r$ to the case $r = \infty$,
too. Namely, in this case consider $n + 1$ horoballs $B_h$ which are mutually tangent. The convex hull of their
base points at infinity will be a totally asymptotic or ideal regular simplex $T_{reg} \in \overline{\mathbf{H}}^n$ of finite
volume. Hence, in this case it is legitimate to write
\begin{equation}
d_n(\infty) = (n + 1)\frac{vol(B_h) \cap T_{reg}}{vol(T_{reg})}. \tag{1.3}
\end{equation}
Then for a horoball packing $\cB_h$, there is an analogue of ball packing, namely (cf. \cite{B78}, Theorem
4)
\begin{equation}
\delta_n(\cB_h, B_h) \le d_n(\infty),~ \forall B_h \in \cB_h. \tag{1.4}
\end{equation}
The upper bound $d_n(\infty)$ $(n=2,3)$ is attained for a regular horoball packing, that is, a
packing by horoballs which are inscribed in the cells of a regular honeycomb of $\overline{\mathbf{H}}^n$. For
dimensions $n = 2$, there is only one such packing. It belongs to the regular tessellation $\{\infty, 3 \}$ . Its dual
$\{3,\infty\}$ is the regular tessellation by ideal triangles all of whose vertices are surrounded
by infinitely many triangles. This packing has in-circle density $d_2(\infty)=\frac{3}{\pi} \approx 0.95493.. $. 

In $\overline{\mathbf{H}}^3$ there is exactly one horoball packing whose Dirichlet--Voronoi cells give rise to a
regular honeycomb described by the Schl\"afli symbol $\{6,~3,~3\}$ . Its
dual $\{3,3,6\}$ consists of ideal regular simplices $T_{reg}$  with dihedral angle $\frac{\pi}{3}$ building up a 6-cycle around each edge 
of the tessellation.

\subsection{Optimal packings by horoballs of different types in $\overline{\mathbf{H}}^n$} 

In \cite{KSz} we have determined the optimal horoball packing arrangements and their densities
for all four fully asymptotic Coxeter tilings (Coxeter honeycombs) in $\overline{\mathbf{H}}^3$. 
Centers of horoballs are required to lie at vertices of the
regular polyhedral cells constituting the tiling. {\it \bf We allow horoballs of
different types at the various vertices}. We have proved that the known B\"or\"oczky--Florian 
density upper bound for "congruent horoball" packings of $\overline{\mathbf{H}}^3$
remains valid for the class of fully asymptotic Coxeter tilings, even if
packing conditions are relaxed by allowing horoballs of different
types under prescribed symmetry groups. The consequences of this remarkable result are discussed for
various Coxeter tilings (see \cite{KSz}), but in this paper we consider only the tetrahedral Coxeter tilings.

The Coxeter tiling $\{ 3,3,6 \}$ is a three-dimensional honeycomb with cells comprised of fully asymptotic
regular tetrahedra. We arbitrarily choose such a tetrahedron $T_{reg}=E_0E_1E_2E_3$,
and place the horoball centers at the ideal vertices $E_0, \dots,
E_3$. We vary the types of the horoballs so that they satisfy our constraints of non-overlap.
The density of horoball packings to the above honeycomb can be computed by the new type of the simplicial density function:
\begin{defn}
The density of a horoball packing in Coxeter honeycomb $\{ 3,3,6 \}$
is defined as
\begin{equation}
\delta(\mathcal{B}_{336})=\frac{\sum_{i=1}^{4} vol(B_i \cap T_{reg})}{vol(T_{reg})}. \tag{1.5}
\end{equation}
\end{defn}
\begin{figure}
\begin{center}
\includegraphics[width=6cm]{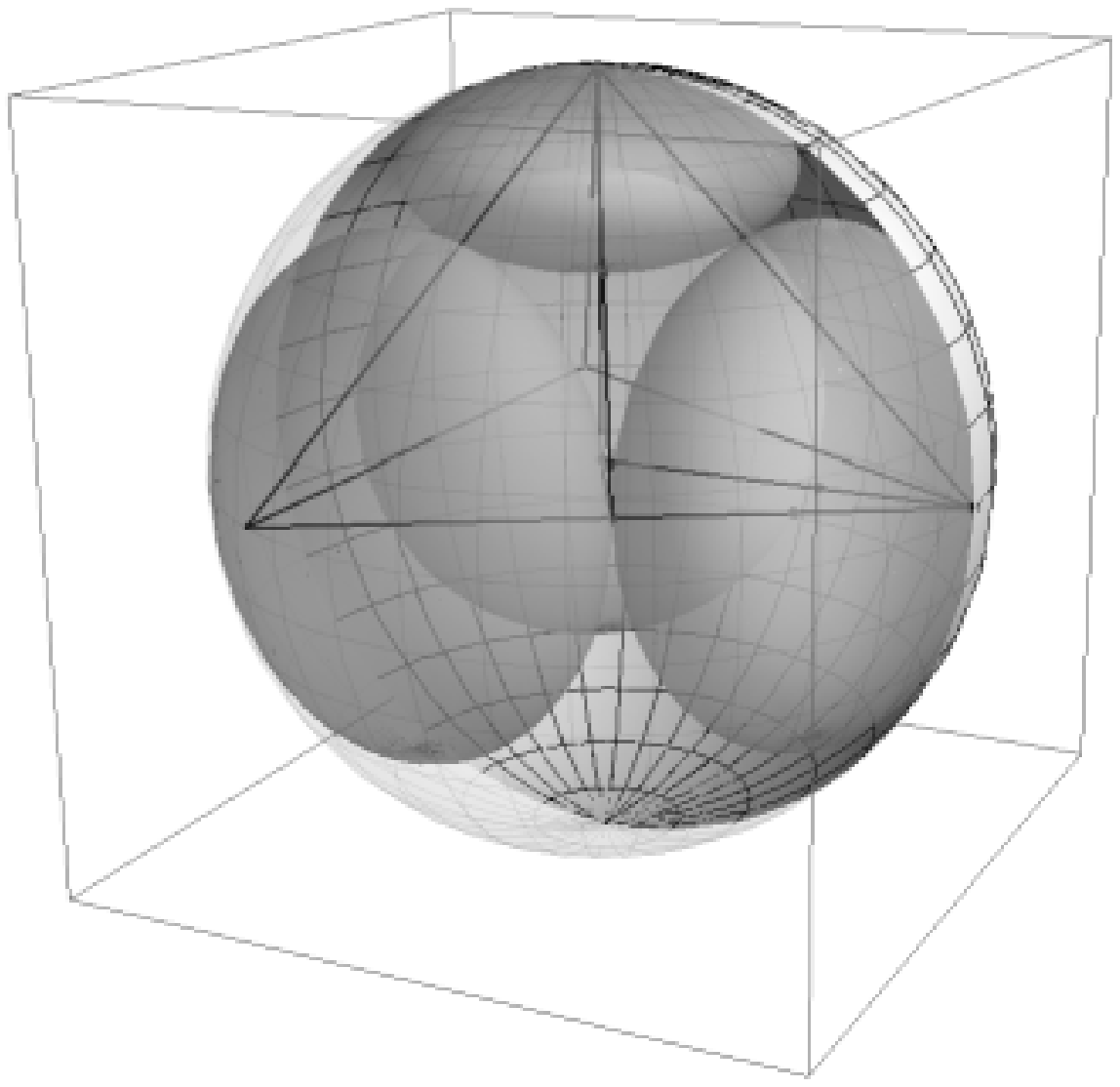}
\includegraphics[width=6cm]{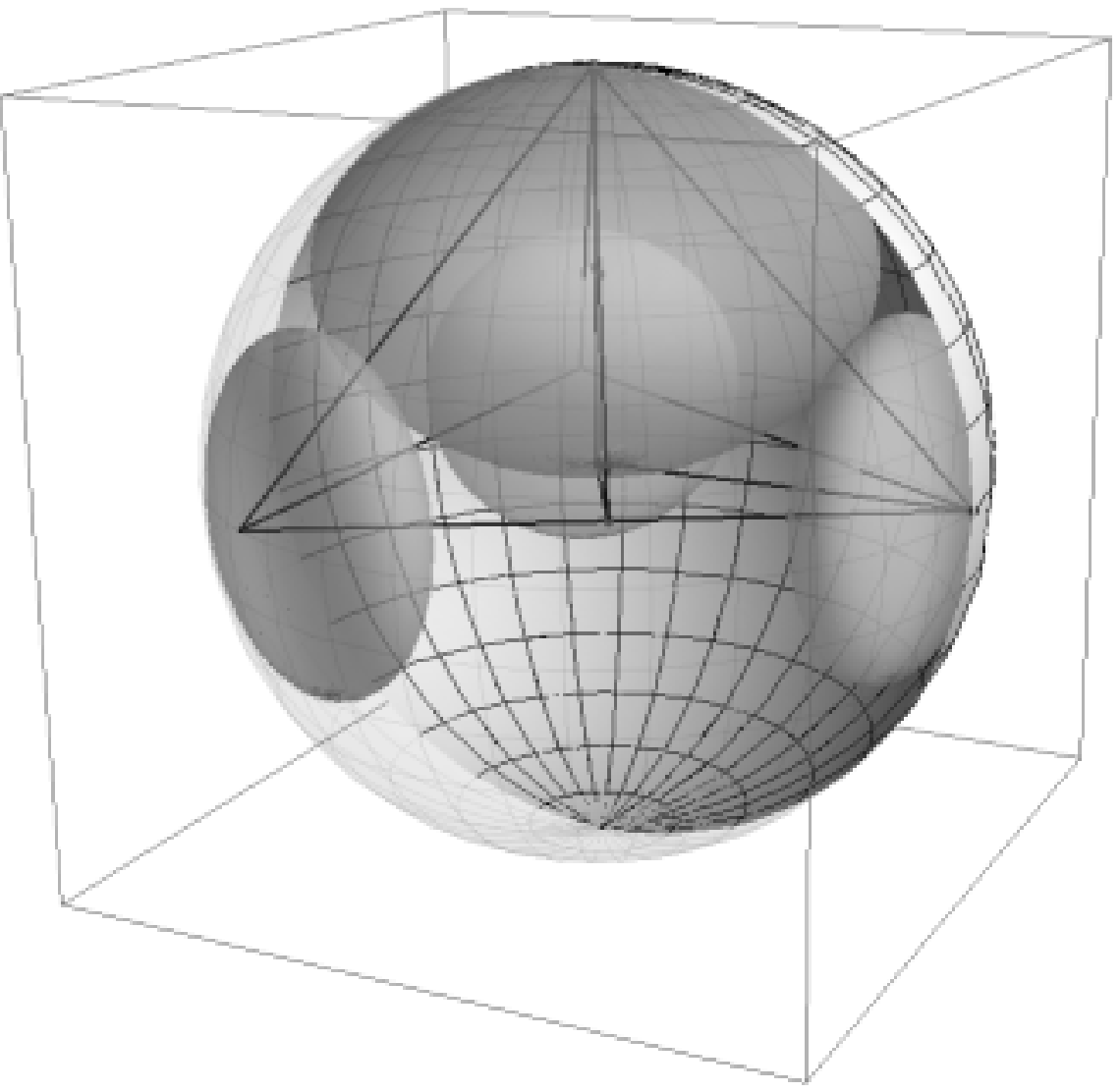}
\\a. ~~~~~~~~~~~~~~~~~~~~~~~~~~~~~~b. \\
\caption{Two optimal horoball arrangemets of $(3,3,6)$ tiling.}
\end{center}
\end{figure}
\begin{theorem}
There are two distinct optimally dense horoball arrangements
$\mathcal{B}_{336}^i,$ $(i=1,2)$ for the tetrahedral Coxeter tiling
$(3,3,6)$ with the same density: $\delta(\mathcal{B}_{336}^i) \approx
0.85327609$.
\end{theorem}
Fig.~1. shows the two optimal horoball arrangements in Beltrami-Cayley-Klein model. 
\subsection{Generalization of the simplicial density function}
\begin{defn}
We consider an arbitrary fully asymptotic simplex $T=E_0 E_1$ $E_2 E_3 \dots E_n$ in the $n$-dimensional hyperbolic space $\overline{\mathbf{H}}^n$.
Centers of horoballs are required to lie at vertices of $T$. We allow horoballs $(B_i, ~ i=1,2,\dots,n)$ of different types at the various vertices and require to form a packing, 
moreover we assume that $$card(B_i \cap [E_{i_0}E_{i_1}\dots E_{i_{n-1}}]) \le 1, ~ ~ i_j \ne i, ~ ~ j \in \{0,1,\dots ,n-1\}.$$
(The hyperplane of points $E_{i_0},$ $E_{i_1},$ $\dots,$ $E_{i_{n-1}}$ is denoted by $[E_{i_0}E_{i_1}\dots E_{i_{n-1}}]$ may touch the horoball $B_{i_n}$.) 
The generalized simplicial density function for the above simplex and horoballs is defined as
\begin{equation}
\delta(\mathcal{B})=\frac{\sum_{i=0}^{n} vol(B_i \cap T)}{vol(T)}. \notag
\end{equation}
\end{defn}
\emph{For $n=3$ the main problem is} to find a 
fully asymptotic tetrahedron $T_{opt} \in\ \overline{\bH}^3$ and horoballs $B_i$ centered at the vertices $E_i$ such that  
the density $\delta(\cB)$ (see Definition 1.3) of the corresponding horoball arrangement 
is maximal. In this case the horoball arrangement $\cB$ is 
said to be \emph{locally optimal.} 

In this paper, we study locally optimal horoball packings for fully asymptotic tetrahedra
in $\overline{\mathbf{H}}^3$, while allowing different
types of horoballs to be centered at the vertices of the fully asymptotic tetrahedra. 

The general investigation of this problem seems to be very difficult, thus 
we assume for simplicity that the tetrahedron has at least one plane of symmetry, simply called plane symmetric. So it will have only one free
angle parameter. Then this ideal tetrahedron has 2 orthogonal symmetry planes and 2 symmetry lines of these planes, i.e. the symmetry group is of order
8, denoted by ${\mathbf{2}~^* \mathbf{2}={\mathbf{\overline{4}m2}}}$ (by Conway and Hermann-Mauguin, respectively).

\section{Computations in projective model}

For $\overline{\mathbf{H}}^n$ $n \geq 2$ we use the projective model in Lorentz space
$\mathbf{E}^{1,n}$ of signature $(1,n)$, i.e.~$\mathbf{E}^{1,n}$ is
the real vector space $\mathbf{V}^{n+1}$ equipped with the bilinear
form of signature $(1,n)$
\begin{equation}
\langle ~ \mathbf{x},~\mathbf{y} \rangle = -x^0y^0+x^1y^1+ \dots + x^n y^n \tag{2.1}
\end{equation}
where the non-zero vectors
$$
\mathbf{x}=(x^0,x^1,\dots,x^n)\in\mathbf{V}^{n+1} \ \  \text{and} \ \ \mathbf{y}=(y^0,y^1,\dots,y^n)\in\mathbf{V}^{n+1},
$$
are determined up to real factors and they represent points in
$\mathcal{P}^n(\mathbf{R})$. ${\mathbf{H}}^n$ is represented as the
interior of the absolute quadratic form
\begin{equation}
Q=\{[\mathbf{x}]\in\mathcal{P}^n | \langle ~ \mathbf{x},~\mathbf{x} \rangle =0 \}=\partial \mathbf{H}^n \tag{2.2}
\end{equation}
in real projective space $\mathcal{P}^n(\mathbf{V}^{n+1},
\mbox{\boldmath$V$}\!_{n+1})$. All proper interior point $\mathbf{x} \in {\mathbf{H}}^n$ are characterized by
$\langle ~ \mathbf{x},~\mathbf{x} \rangle < 0$.

The points on the boundary $\partial \mathbf{H}^n $ in
$\mathcal{P}^n$ represent the absolute points at infinity of $\overline{\mathbf{H}}^n$.
Points $\mathbf{y}$ with $\langle ~ \mathbf{y},~\mathbf{y} \rangle >
0$ lie outside of $\overline{\mathbf{H}}^n$ and are called outer points
of $\mathbf{H}^n $. Let $X([\mathbf{x}]) \in \mathcal{P}^n$ a point;
$[\mathbf{y}] \in \mathcal{P}^n$ is said to be conjugate to
$[\mathbf{x}]$ relative to $Q$ when $\langle ~
\mathbf{x},~\mathbf{y} \rangle =0$. The set of all points conjugate
to $X([\mathbf{x}])$ form a projective polar hyperplane
\begin{equation}
pol(X):=\{[\mathbf{y}]\in\mathcal{P}^n | \langle ~ \mathbf{x},~\mathbf{y} \rangle =0 \}. \tag{2.3}
\end{equation}
Hence the bilinear form $Q$ by (2.1) induces a bijection
(linear polarity $\mathbf{V}^{n+1} \rightarrow
\mbox{\boldmath$V$}\!_{n+1})$
from the points of $\mathcal{P}^n$
onto its hyperplanes.

Point $X [\bold{x}]$ and the hyperplane $\alpha
[\mbox{\boldmath$a$}]$ are called incident if the value of the
linear form $\mbox{\boldmath$a$}$ on the vector $\bold{x}$ is equal
to zero; i.e., $\bold{x}\mbox{\boldmath$a$}=0$ ($\mathbf{x} \in \
\mathbf{V}^{n+1} \setminus \{\mathbf{0}\}, \ \mbox{\boldmath$a$} \in
\mbox{\boldmath$V$}_{n+1} \setminus \{\mbox{\boldmath$0$}\}$).
Straight lines in $\mathcal{P}^n$ are characterized by the
2-subspaces of $\mathbf{V}^{n+1} \ \text{or $(n-1)$-spaces of} \
\mbox{\boldmath$V$}\!_{n+1}$ (see e.g. in \cite{M97}).

In this paper we set the sectional curvature of $\overline{\mathbf{H}}^n$,
$K=-k^2$, to be $k=1$. The distance $s$ of two proper points
$(\mathbf{x})$ and $(\mathbf{y})$ is calculated by the formula:
\begin{equation}
\cosh{{s}}=\frac{-\langle ~ \mathbf{x},~\mathbf{y} \rangle }{\sqrt{\langle ~ \mathbf{x},~\mathbf{x} \rangle
\langle ~ \mathbf{y},~\mathbf{y} \rangle }} . \tag{2.4}
\end{equation}
The foot point $Y(\bold y)$ of the perpendicular, dropped from the point $X(\bold x)$ 
on the plane $(u)$, has the following form:
\[
\bold y=\bold x -\frac{\langle \bold x, \bold u \rangle}
{\langle \bold u, \bold u \rangle} \bold u.
\tag{2.5}
\]
%%%%%%%%%%%%%%%%%%%%%%%%%%%%%%%%%%%%%%%%%%%%%%%%%%%%%%%%%%%%%%%%%

\subsection{On horospheres in hyperbolic space $\overline{\mathbf{H}}^3$}

\begin{defn}
A horosphere in the hyperbolic geometry is the surface orthogonal to
the set of parallel lines, passing through the same point on the absolute quadratic surface 
(simply absolute) of $\overline{\mathbf{H}}^3$.
\end{defn}

We represent $\overline{\mathbf{H}}^3$ by the Beltrami-Cayley-Klein
ball model. We introduce a projective coordinate system using vector
basis $\bold{b}_i \ (i=0,1,2,3)$ for $\mathcal{P}^3$ where the
coordinates of the center of the model is $(1,0,0,0)$. We pick up an
arbitrary point at infinity to be $E_3(1,0,0,1)$ (see Fig.~3).

We obtain the following equation for the horosphere in our
Beltrami-Cayley-Klein model of $\overline{\mathbf{H}}^3$ above:
\begin{equation}
-2 s x^0 x^0-2 x^3 x^3+ 2 (s+1)(x^0 x^3)+(s-1)(x^1 x^1+ x^2 x^2)=0  \tag{2.6}
\end{equation}
\begin{remark}
\begin{enumerate}
\item We get the equation of the horosphere in the Cartesian
coordinate system
($x:=\frac{x^1}{x^0},~y:=\frac{x^2}{x^0},~z:=\frac{x^3}{x^0}$):
\begin{equation}
\frac{2(x^2+y^2)}{1-s}+\frac{4(z-\frac{s+1}{2})^2}{{(1-s)}^2}=1.  \tag{2.7}
\end{equation}
\item We will also use the equation of the horosphere with center $E_2$ (see Fig.~3):
$$\frac{2(\frac{1}{4}x^2+y^2+\frac{3}{4}z^2+\frac{\sqrt{3}}{2}xz)}{1-s_1}+\frac{4(\frac{\sqrt{3}}{2}x-\frac{1}{2}z-\frac{s_1+1}{2})^2}{{(1-s_1)}^2}=1. 
$$
\end{enumerate}
\end{remark}

The length $l(x)$ of a horocycle arc to a chord segment $x$ is determined
by the classical formula due to {{J.~Bolyai}}:
\begin{equation}
l(x)=k \sinh{\frac{x}{k}}.  \tag{2.8}
\end{equation}

The intrinsic geometry of the horosphere is Euclidean, therefore,
the area $\mathcal{A}$ of a horospherical triangle will be computed by
the formula of Heron. The volume of the horoball sectors can be
calculated using another formula by {{J.~Bolyai}}. If the area of a domain on the horoshere is $\mathcal{A}$, 
the volume determined
by $\mathcal{A}$ and the aggregate of axes drawn from $\mathcal{A}$ is equal to
\begin{equation}
V=\frac{1}{2}k\mathcal{A}, \tag{2.9}
\end{equation}
we assume that $k=1$ in this paper.
\subsection{The volume of the fully asymptotic tetrahedron in $\overline{\mathbf{H}}^3$}
A tetrahedron $T$ in hyperbolic space is determined in general by its six dihedral
angles i.e. the mutual angles between the four faces of this tetrahedron.
In order to calculate the volume $vol(T)$ of $T$; one can dissect $T$ into six orthoschemes and their volumes can be calculated by the 
Lobachevsky formula.

A plane orthoscheme is a right-angled triangle, whose area can be expressed by the 
well known defect formula. For three-dimensional spherical orthoschemes, {{Ludwig Schl\"afli},} about 1850,
was able to find the volume differential depending on differentials of the not fixed 
3 dihedral angles.
Already in 1836, {{Lobachevsky}} found 
a volume formula for three-dimensional hyperbolic orthoschemes $\mathcal{O}$ \cite{BH}.

The integration method for orthoschemes of dimension three was generalized by
{{B\"ohm}} in 1962 \cite{BH} to spaces of constant nonvanishing curvature of arbitrary
dimension. 
\begin{figure}[ht]
\centering
\includegraphics[width=10cm]{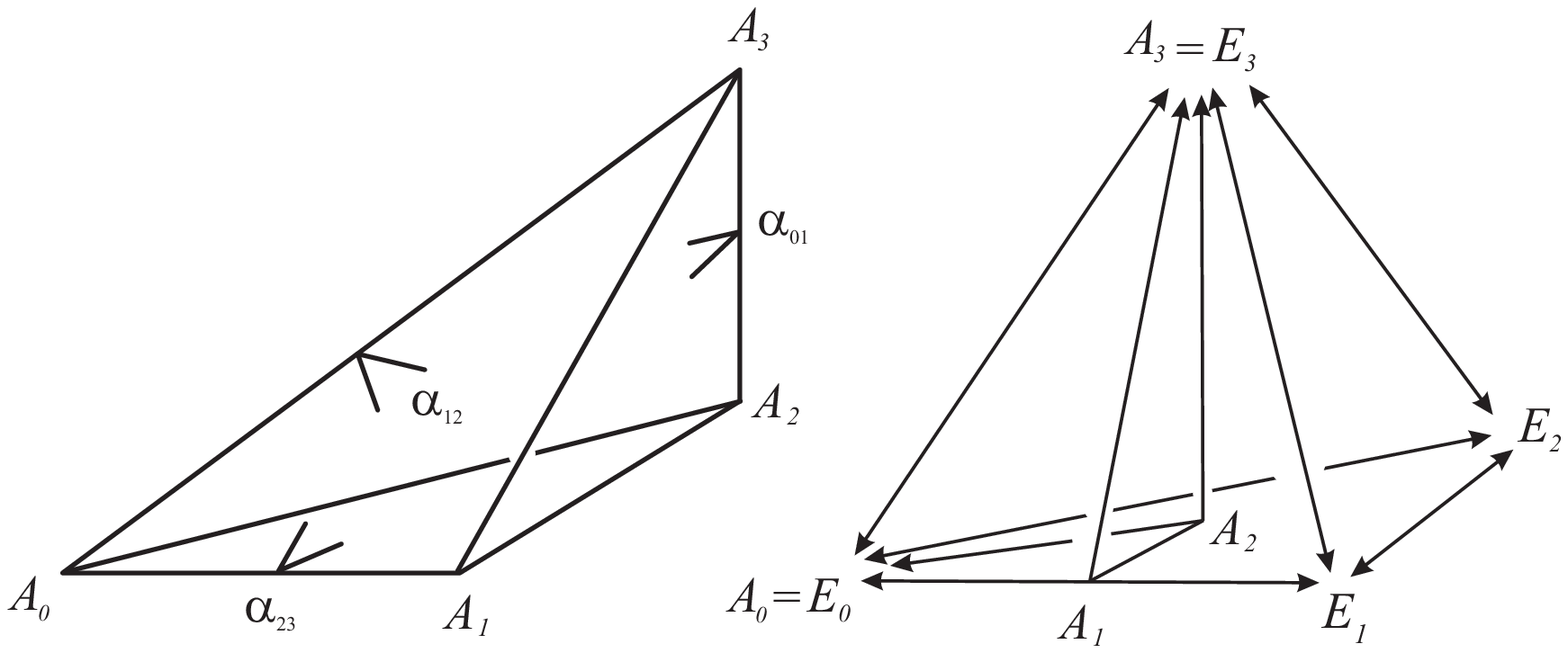}
\caption{}
\label{}
\end{figure}
\begin{theorem}{\rm{(N.~I.~Lobachevsky).}} The volume of a three-dimensional hyperbolic 
ortho\-scheme $\mathcal{O} \subset \overline{\mathbf{H}}^3$ 
is expressed with the dihedral angles $\alpha_{01},\alpha_{12},\alpha_{23}, \ (0 \le \alpha_{ij} \le \frac{\pi}{2})$ 
(Fig.~2) in the following form:

\begin{align}
&vol(\mathcal{O})=\frac{1}{4} \big\{ \mathcal{L}(\alpha_{01}+\theta)-
\mathcal{L}(\alpha_{01}-\theta)+\mathcal{L}\big(\frac{\pi}{2}+\alpha_{12}-\theta\big)+ \notag \\
&+\mathcal{L}\big(\frac{\pi}{2}-\alpha_{12}-\theta\big)+\mathcal{L}(\alpha_{23}+\theta)-
\mathcal{L}(\alpha_{23}-\theta)+2\mathcal{L}\big(\frac{\pi}{2}-\theta\big) \big\}, \notag
\end{align} 
where $\theta \in [0,\frac{\pi}{2})$ is defined by the following formula:
$$
\tan(\theta)=\frac{\sqrt{ \cos^2{\alpha_{12}}-\sin^2{\alpha_{01}} \sin^2{\alpha_{23}
}}} {\cos{\alpha_{01}}\cos{\alpha_{23}}} 
$$
and where $\mathcal{L}(x):=-\int\limits_0^x \log \vert {2\sin{t}} \vert dt$ \ denotes the
Lobachevsky function. 
\end{theorem}

We obtain the volume formula for the asymptotical orthoschems if $A_0$ and $A_3$ are at the absolute i.e. $\alpha:=\alpha_{01}=\frac{\pi}{2}-\alpha_{12}=\alpha_{23}=\theta$:
\[
vol(\mathcal{O_\infty})=\frac{1}{2} \mathcal{L}(\alpha). \tag{2.10}
\]
J.~Milnor's volume formula of a fully asymptotic tetrahedron $T(\alpha,\beta)$ with dihedral angles $\alpha, \beta, \gamma$ is determined by its division 
into orthoschems (see \cite{Mi94}):
\[
vol(T(\alpha,\beta))=\mathcal{L}(\alpha)+\mathcal{L}(\beta)+\mathcal{L}(\gamma),~ ~ \text{where}~ ~ \alpha+\beta+\gamma=\pi. \tag{2.11}
\]
As an easy consequence, we get the following
\begin{lemma}
{The volume formula $vol(T(\alpha))$ of a fully asymptotic plane symmetric tetrahedron $T(\alpha)$ can be derived by the duplication law}:
\[
vol(T(\alpha))=2\mathcal{L}(\alpha)-\mathcal{L}(2\alpha)=-2\mathcal{L}(\alpha+\frac{\pi}{2}). \tag{2.12}
\]
\end{lemma}
In \cite{K98} there are further results for the volumes of orthoschemes
and simplices in higher dimensions. 
\section{Horoball packings for asymptotic tetrahedra}
The aim of this section is to determine the optimal packing arrangement and its
densities for the fully asymptotic tetrahedra in $\overline{\mathbf{H}}^3$.

We will make heavy use of the following Lemma (see also \cite{Sz05-1}):

\begin{lemma}
Let $B_1$ and $B_2$ denote two horoballs with ideal centers $C_1$ and
$C_2$ respectively. Take $\tau_1$ and $\tau_2$ to be two congruent
trihedra, with vertices $C_1$ and $C_2$. Assume that these horoballs
$B_1(x)$ and $B_2(x)$ are tangent at point $I(x)\in {C_1C_2}$ and
${C_1C_2}$ is a common edge of the two trihedra $\tau_1$ and
$\tau_2$. We define the point of contact $I(0)$ such that the
following equality holds for the volumes of horoball sectors:
\begin{equation}
V(0):= 2 vol(B_1(0) \cap \tau_1) = 2 vol(B_2(0) \cap \tau_2). \notag
\label{szirmai-lemma}
\end{equation}
If $x$ denotes the hyperbolic distance between $I(0)$ and $I(x)$,
then the function
\begin{equation}
V(x):= vol(B_1(x) \cap \tau_1) + vol(B_2(x) \cap \tau_2)=\frac{V(0)}{2}(e^{2x}+e^{-2x}) \notag
\end{equation}
strictly increases as~$x\rightarrow\pm\infty$.
\end{lemma}

We arbitrarily choose a fully asymptotic tetrahedron $T(\alpha)=E_0E_1E_2E_3$ with at least one plane of symmetry 
(see Section 1.4 and Fig.~3), and place the horoball centers at vertices $E_0, \dots,
E_3$. We vary the types of the horoballs so that they satisfy our constraints of non-overlap.
The packing density is obtained by Definition 1.3. The dihedral angles of the above tetrahedron at the edges 
$E_0E_1,~E_0E_2,~E_1E_3,~E_2E_3$ are denoted by $\alpha$ and it is clear that the
dihedral angles at the remaining edges $E_0E_3$ and $E_1E_2$ are $\pi-2\alpha$.

We introduce a Euclidean projective coordinate system to the tetrahedron $T(\alpha)$ by the following coordinates of the vertices:

$$E_0:=(1,0,\sqrt{1-z^2},z),~(-1 \le z \le 1); E_1=\big(1,\frac{\sqrt{3}}{2},0,-\frac{1}{2}\big);$$
$$E_2=\big(1,-\frac{\sqrt{3}}{2},0,-\frac{1}{2}\big);~E_3:=(1,0,0,1).$$
$E_0$ lies on the symmetry plane $x=0$ of tetrahedron $T(\alpha)$ (see Fig.~3). 
In order to determine the optimal ball arrangement and its fully asymptotic tetrahedron 
first we consider the following situations:
\begin{enumerate}
\item In an optimally dense packing, the horoball packing must be locally stable, i.e. each ball is fixed by its neighboring horoballs or by the opposite 
face of the tetrahedron. Otherwise the density could be improved by blowing up at least one horoball until it touches a neighboring 
horoball or its opposite face of tetrahedron.
\item We fix the parameter $z$ $(z\in (-1,1)\setminus 0)$ and blow up the horoballs $B_i$ centered at the vertices of $T(\alpha)$ so that the volumes
$vol(B_i \cap T(\alpha)), ~ (i=0,1,2,3)$ are equal until certain two neighbouring horoballs touch each other.
In this situation we distinguish two cases for the varying vertex $E_0$: 

{\bf i}.) If $0 < z < 1$ then the horoballs $B_3-B_0$ and $B_1-B_2$ touch each other the other horoballs do not.

{\bf ii}.) If $-1 < z < 0$ then each of the horoball pairs $B_0-B_1,~B_0-B_2,~B_1-B_3,~B_2-B_3$ have exactly one common point, respectively, and
$B_3-B_0$ and $B_1-B_2$ do not touch each other.

If $z=0$ the tetrahedron is regular $T(\alpha)=T_{reg}, ~ (\alpha=\frac{\pi}{3})$ and each horoball touch all neigboring horoballs. 
This case is studied in details formerly in \cite{KSz}.
\item The volume $vol(T(\alpha)$ is determined by Lemma 2.4, formula (2.12), and by the above machinery of hyperbolic geometry we obtain 
the following formulas:
\begin{equation}
\begin{gathered}
\cos(2\alpha)=-\frac{1+2z}{z-2}, ~ ~ (-1 < z < 1),\\
vol(T(\alpha))=-2\mathcal{L}\Big[\frac{1}{2}\arccos\Big(-\frac{1+2z}{z-2}\Big)+\frac{\pi}{2}\Big]. \tag{3.1}
\end{gathered}
\end{equation}
\end{enumerate}

\vspace{3mm}
{\bf i/1}:~ ($z \in (0,\frac{-2}{13}+\frac{6 \sqrt{3}}{13}]$)
\vspace{3mm}

\begin{enumerate}
\item[a.)] First, we define the tangent point $I(0) \in E_0E_3$ of horoballs $B_0(0)$ and $B_3(0)$ so that 
equalities hold for the volumes of all horoball sectors.

The volume sum of horoball sectors $B_i(0) \cap T(\alpha), ~ (i=0,1,2,3)$ in $T(\alpha)$ is
\begin{equation}
 V(0)=\sum_1^4 (vol(B_i(0) \cap T(\alpha)))=4 vol(B_1(0) \cap T(\alpha)). \notag
\end{equation}
In the introduced coordinate system $B_1(0)$ and $B_2(0)$ touch each other at $(1,0,0,-\frac{1}{2})$.
\begin{figure}[ht]
\centering
\includegraphics[width=7cm]{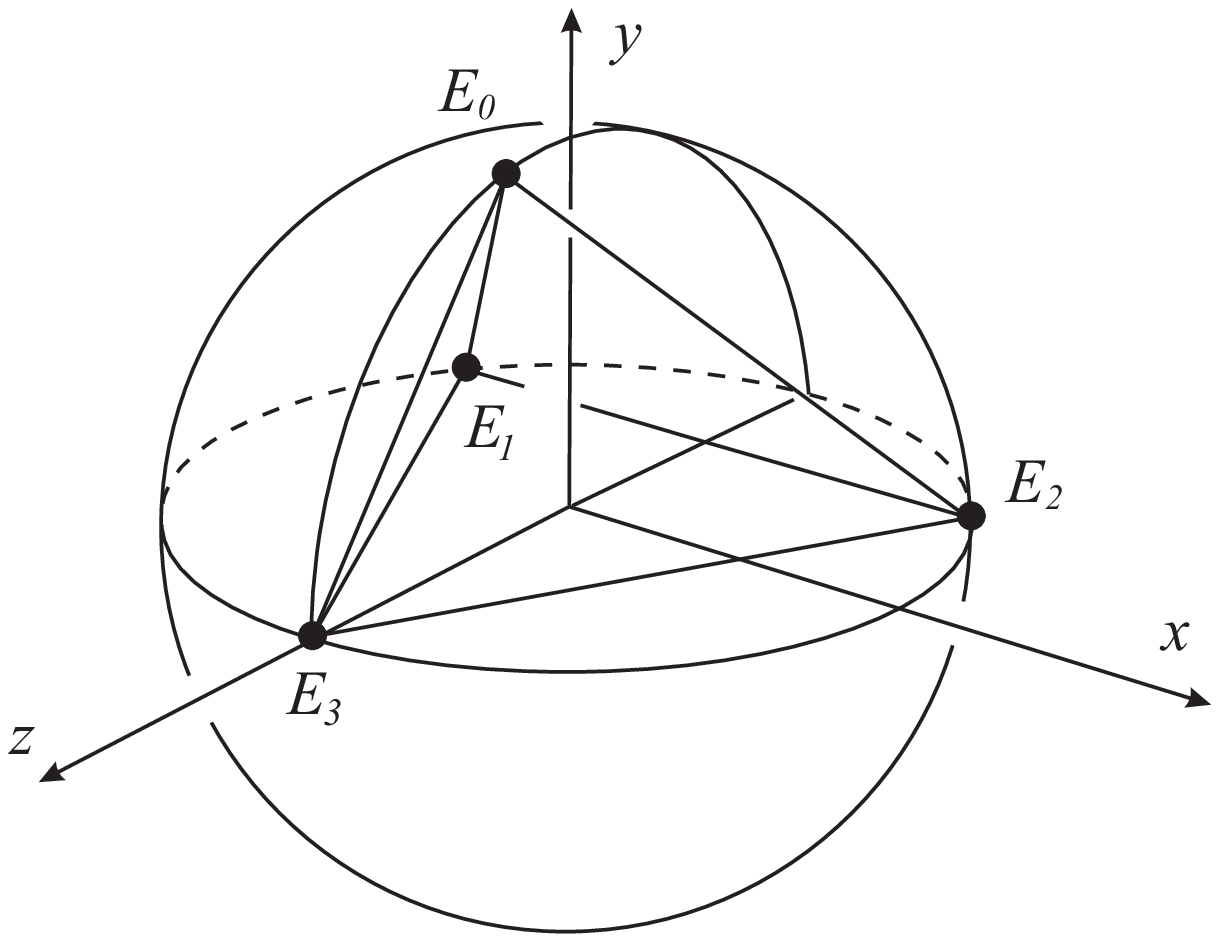}
\caption{}
\label{}
\end{figure}
Consider the point $I(x)$ on the edge $E_0E_3$ where the horoballs $B_i, \ (i=0,3)$ 
are tangent. Here $x$ denotes
the hyperbolic distance between $I(0)$ and $I(x)$. 
Let $I(x_1) \in E_0E_3$ be such a point and parameter $x_1$ where horoball $B_3(x_1)$ touches $B_2(0)$ and $B_1(0)$. 
This horoball arrangement is denoted by $\mathcal{B}_1(x_1(z))$.
Function $V_1(x)$ to $\mathcal{B}_1$ is defined by:
\begin{equation}
V_1(x):= \frac{1}{2}V(0) + \frac{1}{4}V(0) e^{2x}+\frac{1}{4}V(0) e^{-2x}, ~ ~ x \in(0,x_1(z)]. \notag
\end{equation}
By Lemma 3.1 it follows that function $V_1(x)$ strictly increases as
$I(x)$ ($x \in (0,x_1(z)]$) moves away from $I(0)$ along $E_0E_3$.
That means that the density function for any $z \in (0,\frac{-2}{13}+\frac{6 \sqrt{3}}{13}]$,~ 
($\frac{-2}{13}+\frac{6\sqrt{3}}{13} \approx 0.64556191$)

\begin{equation}
\delta(\mathcal{B}_1(x(z)))=\frac{V_1(x)}{vol(T(\alpha))} \notag
\end{equation}
attains its maximum at the endpoint $x_1$. Such a configuration does not occur if $z \in (\frac{-2}{13}+\frac{6\sqrt{3}}{13},1]$.

To study the former density function $\delta(\mathcal{B}_1(x_1(z)))$ we have to compute the volume sum of the horoball sectors 
$V_1(x_1)=\Sigma_1^4(B_i \cap T(\alpha))$ belonging to the parameter 
$x_1$. Of course the volume $V_1(x_1)$ and the volume of the tetrahedron $T(\alpha)$ depend on parameter $z$.
First we calculate the six intersection points of the edges of
$T(\alpha)$ and  horoballs $B_i$ $(i=0,1,2,3)$, using the above projective coordinate system and 
the equations of the horospheres derived from (2.7). These intersection points on the edge $E_iE_j$ is denoted by $M^1_{ij}$ $(i < j), ~
i,j \in \{0,1,2,3\}$. 
\begin{equation}
\begin{gathered}
M^1_{03}=\Big(1,0,-3\frac{\sqrt{1-z^2}}{2z-5},-\frac{z+2}{2z-5}\Big),~ ~ M^1_{13}=(1,-\frac{\sqrt{3}}{4},0,\frac{1}{4}\Big),\\
M^1_{23}=\Big(1,\frac{\sqrt{3}}{4},0,\frac{1}{4}\Big), ~ ~ M^1_{12}=\Big(1,0,0,-\frac{1}{2}\Big),\\
M^1_{02}=\Big(1,\frac{\sqrt{3}(z+2)}{2(z+5)},3\frac{\sqrt{1-z^2}}{z+5},\frac{5z-2}{2(z+5)}\Big), \\
M^1_{01}=\Big(1,-\frac{\sqrt{3}(z+2)}{2(z+5)},3\frac{\sqrt{1-z^2}}{z+5},\frac{5z-2}{2(z+5)}\Big).
\end{gathered} \tag{3.2}
\end{equation}
The volumes of horoball sectors (depending on parameter $z$) can be computed by Bolyai's formulas (2.8-9) and the volume 
of the tetrahedron $T(\alpha)$ is determined by formula (3.1).  
Finally, we obtain the density function  $\delta(\mathcal{B}_1(x_1(z)))=\frac{V_1(x_1)}{vol(T(\alpha))}$ which depends only on parameter 
$z \in (0,\frac{-2}{13}+\frac{6 \sqrt{3}}{13}]$. By careful computer analysis of the above density function we get that the 
function is convex, it attains its maximum at the upper endpoint of the interval $(0,\frac{-2}{13}+\frac{6 \sqrt{3}}{13}]$. 
We have studied the above function with {\it Maple} using that the conditions of the Lobachevsky function are well known \cite{Z07}.    
The graph of $\delta(\mathcal{B}_1(x_1(z)))$ can be seen
in Fig.~4.
\begin{figure}[ht]
\centering
\includegraphics[width=7cm]{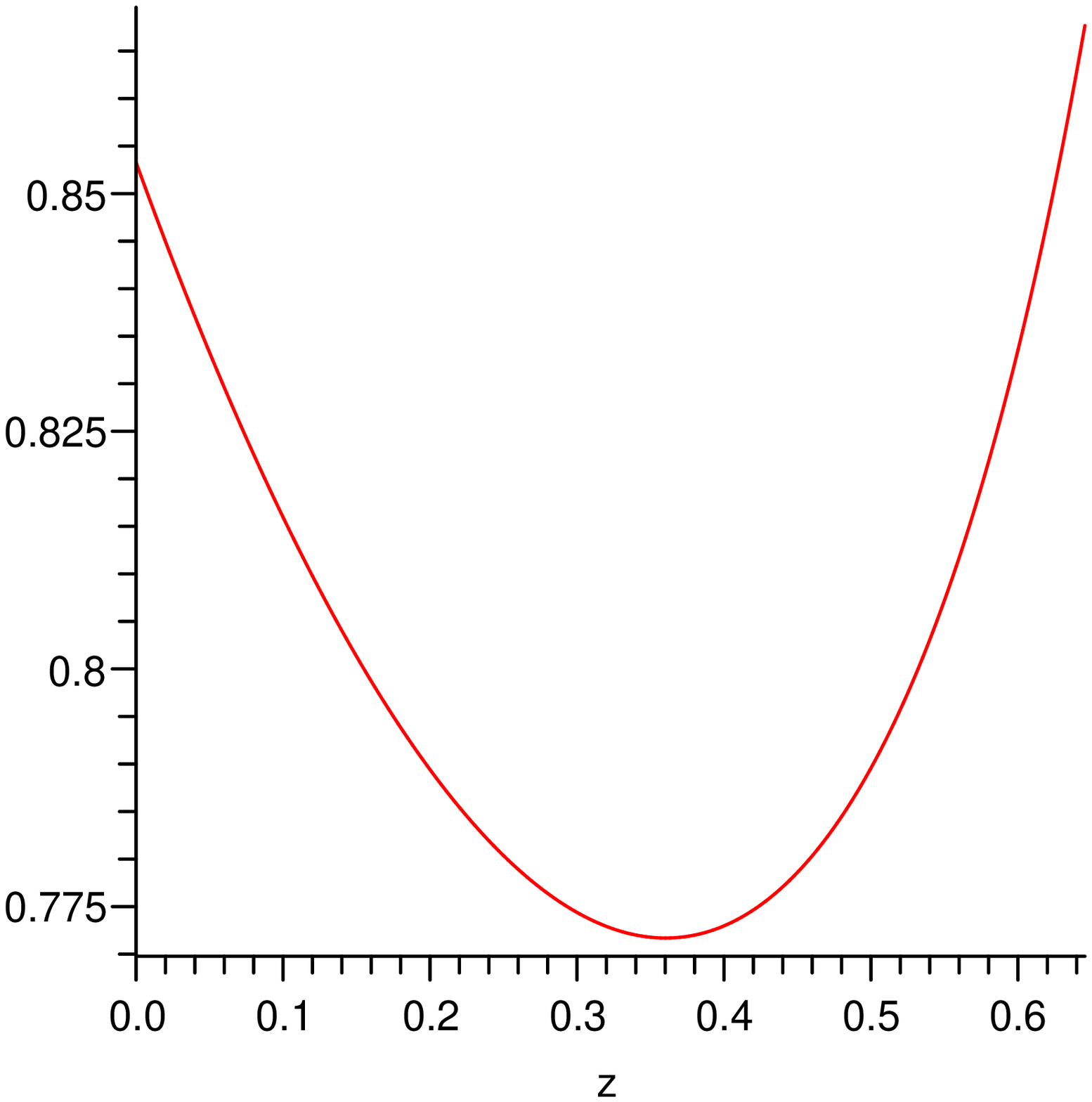}
\caption{}
\label{}
\end{figure}
Note, that the case $z=0$ belongs to the optimal arrangement of B\"or\"oczky--Florian, its density is $\approx 0.85327609$.
The density at the upper endpoint of the interval $(0,\frac{-2}{13}+\frac{6 \sqrt{3}}{13}]$ is larger than that, namely:
\begin{equation}
\delta\Big(\mathcal{B}_1\Big(x_1\Big(\frac{-2}{13}+\frac{6 \sqrt{3}}{13}\Big)\Big)\Big) \approx 0.86767481, \ \ (\alpha \approx 1.30899694). \tag{3.3}
\end{equation}
\item[b.)] We consider the former horoball arrangement $\mathcal{B}_1(x_1(z))$ for any  
$z \in (0,\frac{-2}{13}+\frac{6 \sqrt{3}}{13}]$. We expand the above "larger horoball" $B_3(x_1)$ until it touches the opposite face $E_0E_1E_2$
of the tetrahedron while keeping the other horoballs tangent to it. The arising horoball arrangement are realizable for any 
$z \in (0,\frac{-2}{13}+\frac{6 \sqrt{3}}{13}].$

Similarly to the case a.) we consider a point $I(x)$ $(x \ge x_1)$ (on the edge $E_0E_3$ where the horoballs $B_i, \ (i=0,3)$ 
are tangent at point $I(x) \in E_0E_3$. 
Let $x-x_1$ denote the hyperbolic distance between $I(x_1)$ and $I(x)$. 
Furthermore, let $I(x_2) \in E_0E_3$ be such a point where horoball $B_3(x_2)$ of parameter $x_2$ touches the face $E_0E_1E_2$. 
The "largest horoball" determines the configuration of all other horoballs and this horoball arrangement is denoted by $\mathcal{B}_2(x_2(z))$.
Function $V_2(x)$ is defined as follows:
\begin{equation}
V_2(x):= \frac{1}{2}V(0)e^{-2(x-x_1)} + \frac{1}{4}V(0) e^{2x}+\frac{1}{4}V(0) e^{-2x}, ~ ~ x \in[x_1,x_2]. \notag
\end{equation}
The following Lemma is obtained by examining the above function:
\begin{lemma}
The maxima of function $V_2(x)$ are realized to parameters $x_1$ or $x_2$.
\end{lemma}
By the above Lemma 3.2 it follows that it is sufficient to consider the volume function $V_2(x)$ and the 
density function $\delta(\mathcal{B}_2(x(z)))$ at the parameters $x_1$ and $x_2$. The case $x_1$ was taken in a.) thus we have to consider 
only the case $x=x_2$.

The horoball $B_3(x_2)$ has to touch the side face $E_0E_1E_2$ of fully asymptotic tetrahedron $T(\alpha)$. Thus, it is passing trough 
the  foot point $E_3'$ perpendicularly dropped from the point $E_3$ 
on the plane $E_0E_1E_2$,
$$
E_3'=\Big(1,0,\frac{3\sqrt{(1-z^2)}(1+2z)}{(z+5)^2},-\frac{7z^2+4z-2}{2z^2-4z-7}\Big). 
$$

To examine the density function $\delta(\mathcal{B}_2(x_2(z)))$ we have to compute the sum of the volumes of the horoball sectors 
$V_2(x_2)=\Sigma_1^4(B_i \cap T(\alpha))$ to the parameter 
$x_2$. The volume $V_2(x_2)$ and the volume of the tetrahedron $T(\alpha)$ depend on parameter $z$.
Similarly to the case a.) we have to determine the seven intersection points of the edges of
$T(\alpha)$ and  horoballs $B_i$ $(i=0,1,2,3)$. The intersection points on the edges $E_0 E_1$, $E_0 E_2$, $E_0 E_3$, $E_3 E_2$ and $E_3 E_1$ are 
denoted by $M^2_{ij}$ $(i < j), ~i,j \in \{0,1,2,3\}$. 
On the edge $E_1E_2$ there are no touching point but there are two point of intersections $M^2_{{12}_{1,2}}$.
\begin{equation}
\begin{gathered}
M^2_{03}=\Big(1,0,18\frac{\sqrt{1-z^2}(1+z)}{4z^2+22z+22},\frac{19z^2+22z+4}{4z^2+22z+22}\Big),\\ 
M^2_{13}=\Big(1,-\frac{6\sqrt{3}(z^2-1)}{11z^2-4z-16},0,-\frac{7z^2+4z-2}{11z^2-4z-16}\Big),\\
M^2_{23}=\Big(1,\frac{6\sqrt{3}(z^2-1)}{11z^2-4z-16},0,-\frac{7z^2+4z-2}{11z^2-4z-16}\Big),\\ 
M^2_{{12}_{1,2}}=\Big(1,\pm \frac{\sqrt{3}(13z^2+4z-8)}{2(11z^2-4z-16)},0,-\frac{1}{2}\Big),\\
\end{gathered} \tag{3.4}
\end{equation}
\begin{equation}
\begin{gathered}
M^2_{02}=\Big(1,\frac{2\sqrt{3}(z^2-1)}{4z^2-z-6},-\frac{\sqrt{1-z^2}(z+2)}{4z^2-z-6},-\frac{3z^2+2z-2}{4z^2-z-6}\Big), \\
M^2_{01}=\Big(1,-\frac{2\sqrt{3}(z^2-1)}{4z^2-z-6},-\frac{\sqrt{1-z^2}(z+2)}{4z^2-z-6},-\frac{3z^2+2z-2}{4z^2-z-6}\Big).
\end{gathered} \notag
\end{equation}
The volumes of horoball sectors (depending on parameter $z$) can be computed by Bolyai's formulas (2.8-9) and the volume 
of the tetrahedron $T(\alpha)$ is determined by formula (3.1). 
We get the density function  $\delta(\mathcal{B}_2(x_2(z)))=\frac{V_2(x_2)}{vol(T(\alpha))}$ which depend only on parameter 
$z \in (0,\frac{-2}{13}+\frac{6 \sqrt{3}}{13}]$. By careful investigation (see \cite{Z07}) of the above density function we get that the function 
attains its maximum at the point $z=\frac{-2}{13}+\frac{6 \sqrt{3}}{13}$. The graph of $\delta(\mathcal{B}_2(x_2(z)))$ is shown in Fig.~5.
\begin{figure}[ht]
\centering
\includegraphics[width=7cm]{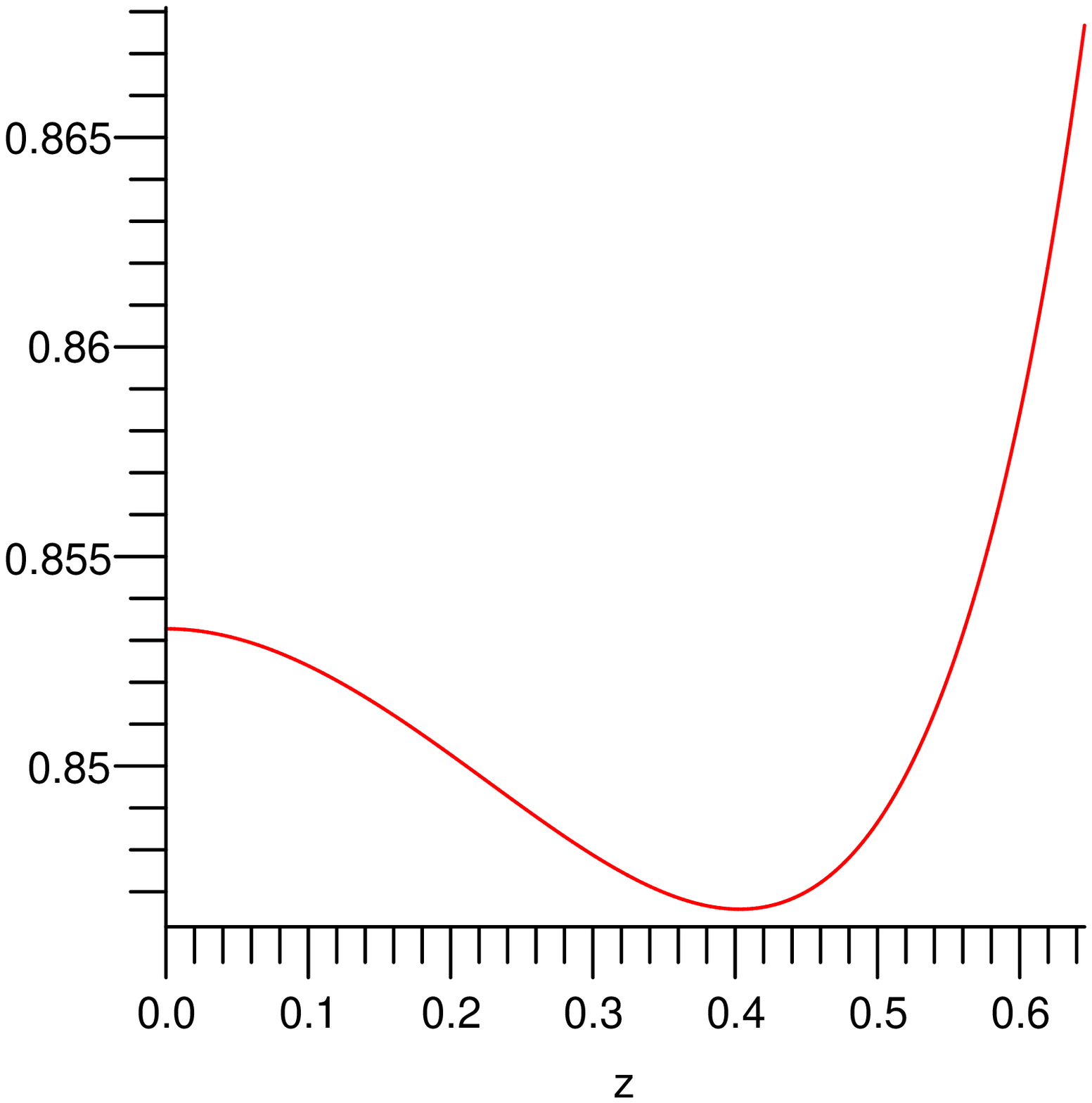}
\caption{}
\label{}
\end{figure}
The density at the point $\frac{-2}{13}+\frac{6 \sqrt{3}}{13}$ is larger than the B\"or\"oczky-Florian upper density, namely
\begin{equation}
\delta\Big(\mathcal{B}_2\Big(x_2\Big(\frac{-2}{13}+\frac{6 \sqrt{3}}{13}\Big)\Big)\Big) \approx 0.86767481, \ \ (\alpha \approx 1.30899694). \tag{3.5}
\end{equation}
The horoball density of of B\"or\"oczky--Florian belongs to the parameter $z=0$. To the parameter $z=\frac{-2}{13}+\frac{6 \sqrt{3}}{13}$
we get the same horoball density as in case a.).
\end{enumerate}
\begin{figure}[ht]
\centering
\includegraphics[width=7cm]{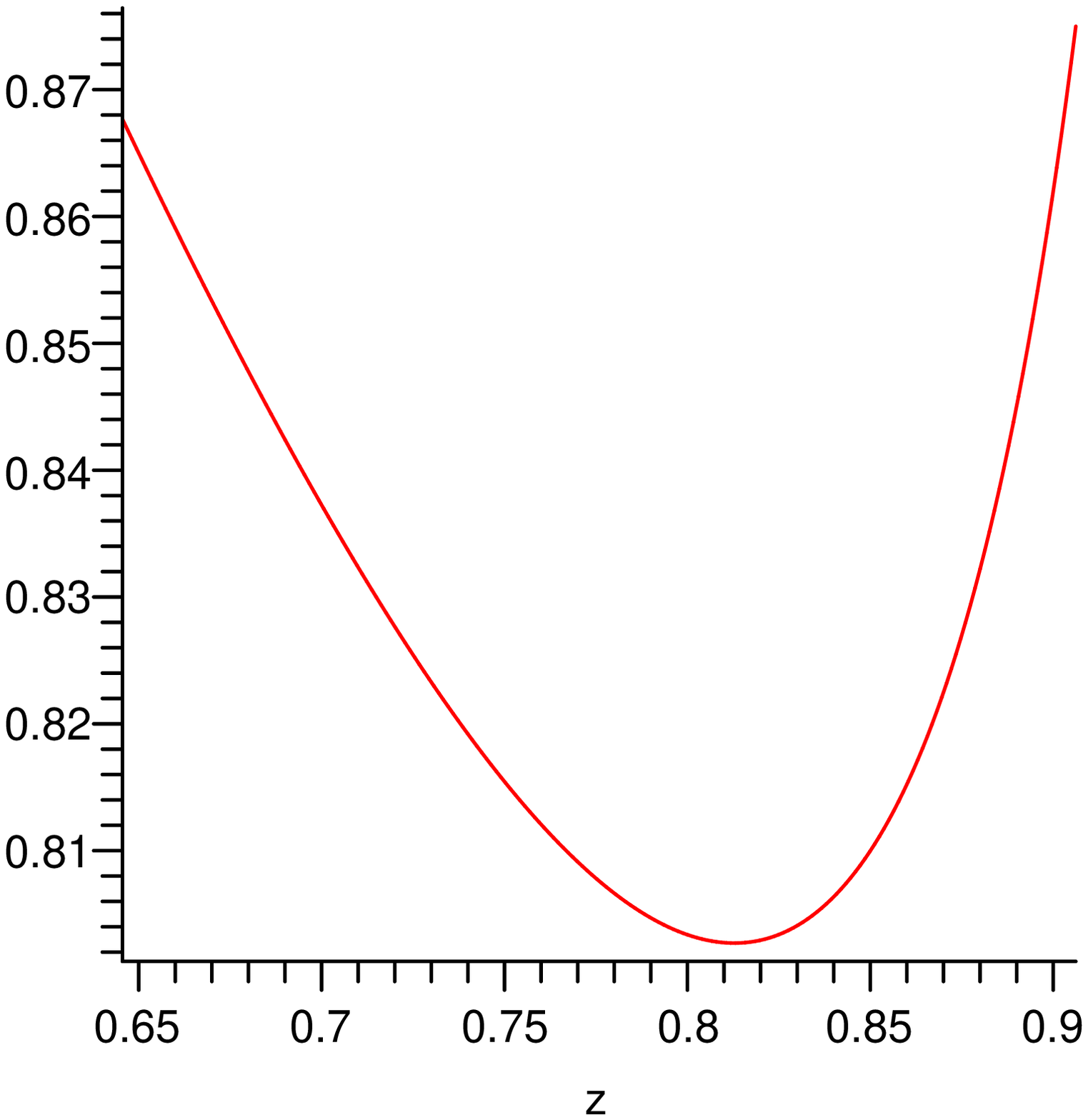}
\caption{}
\label{}
\end{figure}
\vspace{3mm}
{\bf i/2}:~ ($z \in (\frac{-2}{13}+\frac{6 \sqrt{3}}{13},1]$)
\vspace{3mm}

\begin{enumerate}
\item[a.)]  $z \in \Big( \frac{-2}{13}+\frac{6 \sqrt{3}}{13},z_3\Big]$, where
$z_3:= \frac{4}{73} \sqrt{3}\sqrt{178} \cos(\frac{1}{3} \arctan (\frac{73}{8001} 
\sqrt{3} \sqrt{229}))-\frac{26}{73} \approx 0.9061774494.$

In the former case {\bf i/1/b} we have already considered the ball arrangement $\mathcal{B}_1(x_1(z))$ for any  
$z \in (0,\frac{-2}{13}+\frac{6 \sqrt{3}}{13}]$. There we have expanded the "larger horoball" $B_3(x_1)$ until it touches the opposite face $E_0E_1E_2$
of the tetrahedron while keeping other horoballs tangent to it. 
But if $z \in \big(\frac{-2}{13}+\frac{6 \sqrt{3}}{13},z_3\big]$ then we cannot apply this procedure. Namely,
if the horoball $B_3(x_2)$ touches the opposite side face, then $B_3(x_2)$ does not touch $B_1(0)$ and $B_2(0)$.
Then we expand the horoball $B_2$ until it touches $B_3(x_2)$ or the opposite side face.
This means that the "larger horoball" $B_3(x_2)$ touches 
the opposite face $E_0E_1E_2$ at point $E_3'$ and the horoballs $B_2$, $B_0$ touch $B_3$. Then $B_1$ touches only the horoball $B_2$. 
This configuration occurs until the horoball $B_2$ is not tangent to its opposite face $E_3E_1E_0$. 

This arrangement exactly exists if $z \in \big(\frac{-2}{13}+\frac{6 \sqrt{3}}{13},z_3\big]$. 
The horoball $B_2$ touches the opposite side face $E_0E_1E_3$ of $T(\alpha)$ at the point
$E_2':$
\begin{equation}
E_2'=\Big(\frac{5z+7}{2(z+2)},-\frac{\sqrt{3}(2z+1)}{2(z+2)},\frac{3\sqrt{1-z^2}}{2(z+2)},\frac{2z+1}{2(z+2)}\Big). \tag{3.6}
\end{equation}
Similarly to the case a.) the "largest horoball" determines the configuration of all other horoballs and this horoball arrangement is denoted by 
$\mathcal{B}_3(x_3(z))$. The function $V_3(x_3(z))$ and related density function 
$$\delta(\mathcal{B}_3(x_3(z)))=\frac{V_3(x_3(z))}{vol(T(\alpha))}$$ can be analysed in the same way as in {\bf i/1/a} and {\bf i/1/b} and we obtain 
that the convex density function attains its maximum at the above $z_3$.

The graph of $\delta(\mathcal{B}_3(x_3,z)), z \in [\frac{-2}{13}+\frac{6 \sqrt{3}}{13}),z_3]$ is shown in the Fig.~6. While Fig.~7 shows the 
larger horoballs $B_3,~B_2 \in \mathcal{B}_3(x_3(z_3))$ with tetrahedron $T(\alpha)$ in the Beltrami-Cayley-Klein model. 
\begin{figure}[ht]
\centering
\includegraphics[width=8cm]{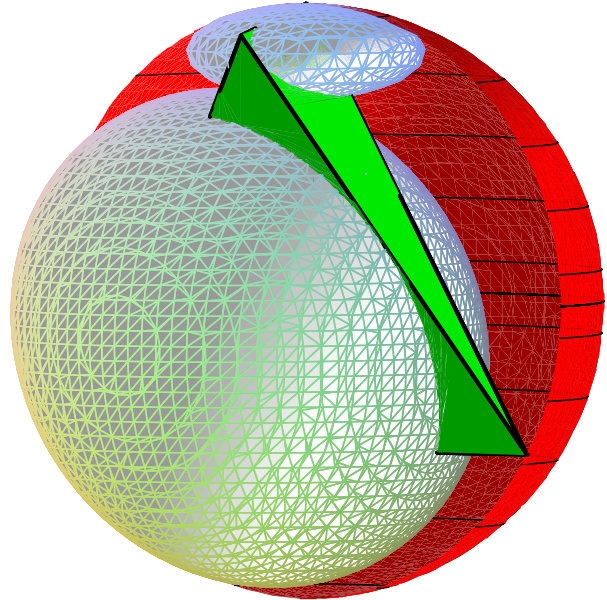}
\caption{The optimal horoball packing arrangement related to a fully asymptotic tetrahedron}
\label{}
\end{figure}
The density at the point $z=z_3$ is larger than the above maxima
\begin{equation}
\delta(\mathcal{B}_3(x_3(z_3)))\approx 0.87499429, \ \ (\alpha \approx 1.44340117). \tag{3.7}
\end{equation}

\item[b.)] $z \in (z_3,1)$

We consider such  ball arrangements $\mathcal{B}_4(x_4(z))$ where the "larger horoballs" $B_3$ and $B_2$ touch their opposite faces $E_0E_1E_2$ and
$E_0E_1E_3$, respectively. The horoballs $B_0,B_3$ and $B_1,B_2$ pairwise touch each other. This horoball packing can be studied similarly to the above cases 
and we can analyse the density function $\delta(\mathcal{B}_4(x_4(z)))$ (see Fig.~8), not detailed further. 
We obtain, that the density function is stricly decreasing in the interval $(z_3,1)$. Thus, there is no horoball arrangement  
with larger density than $\mathcal{B}_3(x_3(z_3))$
with density $\approx 0.87499429$ (see (3.7)).
\begin{figure}[ht]
\centering
\includegraphics[width=7cm]{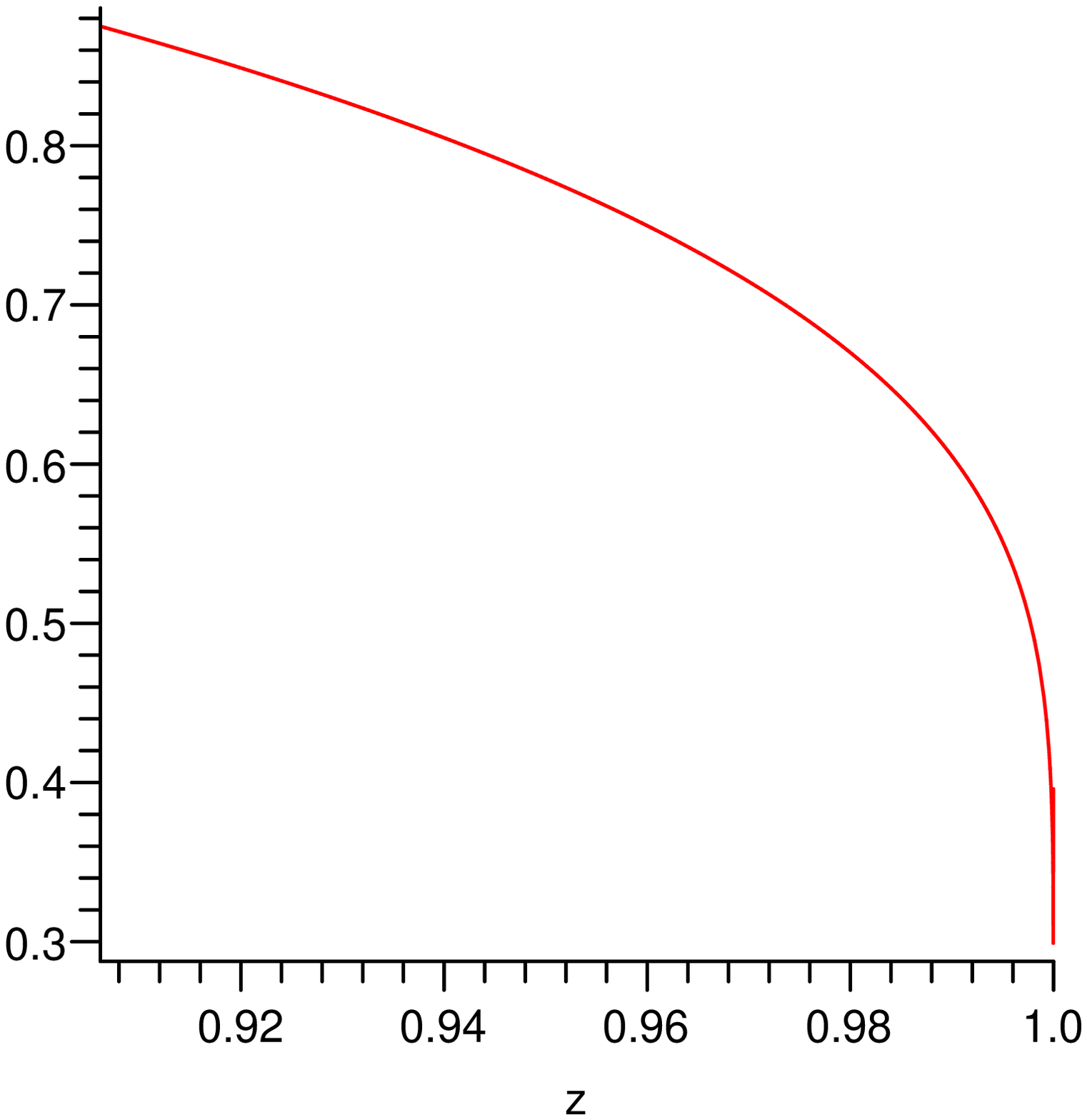}
\caption{}
\label{}
\end{figure}
\end{enumerate}
Analogously to ({\bf i}) we have investigated the cases ({\bf ii}) $(-1 < z < 0)$, not detailed more, and we have obtained the following
\begin{theorem}
The densest horoball arrangement, related to the generalized simplicial density function (Definition 1.3) belongs to the
horoball arrangement $\mathcal{B}_3(x_3(z_3))$ (see (3.7) and Fig.~7) with density $$\delta(\mathcal{B}_3(x_3(z_3)))\approx 0.87499429.$$
\end{theorem}
\begin{rmrk}
This packing seems not to have global extension to the entire hyperbolic space $\overline{\mathbf{H}}^3$, so that 
the same density occurs in each asymptotic tetrahedron. We plan to investigate this problem with E. Moln\'ar and I. Prok by \cite{MPSz06}.
\end{rmrk}

Optimal sphere packings in other homogeneous Thurston geometries
represent a class of open mathematical problems. For
these non-Euclidean geometries only very few results are known
\cite{Sz07-2}, \cite{Sz10-1}. Detailed studies are the objective of
ongoing research.

%%%%%%%%%%%%%%%%%%%%%%%%%%%%%%%%%%%%%%%%%%%%%%%%%%%%%%%%%%%%%%%%%%%%%%%%%%%%%%%%%%%%%%%%%%%%%

%============================================================================%
%                                references                                  %
%============================================================================%

\end{document}